\documentclass[a4paper,twoside,10pt,leqno]{article}
\usepackage{amsmath,amsthm,amssymb,enumerate}
\usepackage{epsfig, float, afterpage, array}

\swapnumbers
\theoremstyle{plain}
\newtheorem{theorem}{Theorem}[section]
\newtheorem{lemma}[theorem]{Lemma}
\newtheorem{corollary}[theorem]{Corollary}

\newtheorem{Forester}[theorem]{The sliding lemma}
\newtheorem{compression}[theorem]{The compressing lemma}

\newtheorem{retraction}[theorem]{The retraction lemma}
\newtheorem{almost}[theorem]{The almost stability theorem \rm \cite[Theorem~III.8.5]{DicksDunwoody89}}

\theoremstyle{definition}
\newtheorem{definition}[theorem]{Definition}
\newtheorem{definitions}[theorem]{Definitions}
\newtheorem{example}[theorem]{Example}

\newtheorem{conventions}[theorem]{Conventions}
\newtheorem{remark}[theorem]{Remark}

\newtheorem{hypotheses}[theorem]{Hypotheses}

\numberwithin{equation}{theorem}

\newcommand{\gen}[1]{\langle#1\rangle}

\DeclareMathOperator{\degg}{height}
\DeclareMathOperator{\length}{length}
\DeclareMathOperator{\ad}{ad}

\DeclareMathOperator{\core}{core}
\def \naturals{\mathbb{N}}
\def \integers {\mathbb{Z}}
\def \reals{\mathbb{R}}

\def\d1{\discretionary{-}{}{-}}

\begin{document}

\pagestyle{myheadings}
\markboth{Retracts of trees and the almost stability theorem}{Warren Dicks and M.\ J.\ Dunwoody}
\title{Retracts of vertex sets of trees \\ and  the  almost stability theorem}

\author{Warren Dicks and M.\ J.\ Dunwoody}

\date{\today}

\maketitle

\begin{abstract}  Let $G$ be a group, let $T$ be an (oriented) $G$-tree with finite
edge stabilizers, and let $VT$ denote the vertex set of $T$.  We show that, 
for each $G$-retract $V'$ of the $G$-set~$VT$,  there exists a $G$-tree 
whose edge stabilizers are finite and whose vertex set is $V'$.  This fact leads to 
various new consequences of the almost stability  theorem.

We also give an example of a group $G$, a $G$-tree $T$ and a $G$-retract $V'$ of $VT$
 such that no $G$-tree has vertex set $V'$.

\medskip

{\footnotesize

}

\end{abstract}

\section{Outline}
\label{sec:intro}

Throughout the article, let $G$ be a group, and let $\naturals$ 
denote the set of finite cardinals, $\{0,1,2,\ldots\}$.
All our $G$-actions will be on the left.

The following extends Definitions~II.1.1 of~\cite{DicksDunwoody89} 
(where $A$ is assumed to have trivial $G$-action).

\begin{definition}\label{def:almost} Let $E$ and $A$ be $G$-sets.

Let $(E,A)$ denote the set of all functions from $E$ to $A$.
An element $v$ of $(E,A)$ has the form $v \colon E \to A$, $e \mapsto v(e)$. 
There is a natural $G$-action on $(E,A)$ such that $(gv)(e) := g(v(g^{-1}e))$
for all $v \in (E,A)$, $g \in G$, $e \in E$. 

Two elements $v$ and $w$ of $(E,A)$ are 
said to be {\it almost equal} if the set 
$$\{e \in E \mid v(e) \ne w(e)\}$$ is finite.
Almost equality is an equivalence relation; the equivalence
classes are called the {\it almost equality classes in $(E,A)$}.

A subset $V$ of $(E,A)$  is said to 
be {\it $G$-stable} if $V$ is closed under the $G$-action.
In general, a $G$-stable subset is the same as a $G$-subset.
\hfill\qed
\end{definition}

In this article, we wish to strengthen the following result.

\begin{almost}\label{thm:almost}  If $E$ is a $G$-set with finite stabilizers, 
and $A$ is a nonempty set with trivial $G$-action, and $V$ is a $G$-stable almost equality class in the 
$G$-set $(E,A)$, then there exists a $G$-tree with finite edge stabilizers 
and vertex set $V$.\hfill\qed
\end{almost}

In the light of Bass-Serre theory, 
the almost stability theorem can be thought of as a broad generalization of 
Stallings' ends theorem.

Let us now recall the notion of a $G$-retract of a $G$-set.
The following alters Definition~III.1.1 of~\cite{DicksDunwoody89} slightly.

\begin{definition}  
A {\it $G$-retract} $U$ of a $G$-set $V$ is a $G$-subset  of  $V$ with the property that,
for each $w \in V - U$, there exists $u \in U$ such that $G_w \le G_u$, or,
equivalently, with the property that there exists a $G$-map,  
called a {\it $G$-retraction}, from $V$ to $U$ which is the identity on $U$.   \hfill\qed
\end{definition}

Chapter IV of~\cite{DicksDunwoody89} collects together a wide variety of consequences
of the almost stability theorem~\ref{thm:almost}.  In some of these applications,
the conclusions assert that certain naturally arising $G$-sets are $G$-retracts of 
vertex sets of  $G$-trees with finite edge stabilizers.  This leads to the question
of whether or not the class of vertex sets of  $G$-trees with finite edge stabilizers is closed under
taking $G$-retracts.   We are now able to answer this in the affirmative;
in Section~\ref{sec:retract} below, we prove that any $G$-retract of the 
 vertex set of a $G$-tree with finite edge stabilizers
is itself the vertex set of a $G$-tree 
with finite edge stabilizers.

In Section~\ref{sec:applications},
we record  the resulting generalizations of the 
almost stability theorem and the applications which are affected.  
In the most classic example,
if $G$ has cohomological  dimension one, 
and $\omega \integers G$ is the augmentation ideal of the group
ring $\integers G$,
one can deduce that $G$ acts freely on a tree whose vertex set is the $G$-set 
$1+\omega\integers G$, and, hence, 
$G$ is a free group; this is a slightly more detailed version of a theorem of 
Stallings and Swan. 

In Section~\ref{sec:generalized}, we record an even more general form of the 
almost stability theorem in which the $G$-action on $A$ need not be trivial.

In Section~\ref{sec:example}, we construct
 a group $G$ and a $G$-retract of a vertex set of a $G$-tree (with infinite edge stabilizers) that is
not itself the vertex set of a $G$-tree.

\section{Operations on trees}

Throughout this section we will be working with the following.

\begin{hypotheses}\label{hyp:sec2}
Let $T = (T,V,E, \iota,\tau)$ be a $G$-tree,
as in~\cite[Definition~I.2.3]{DicksDunwoody89}.  

We write $VT = V$ and $ET = E$, and we view the underlying $G$-set of $T$ as the
disjoint union of $V$ and $E$, written $T = V \vee E$.  Here $\iota\colon E \to V$ is the
{\it initial vertex} map and $\tau\colon E \to V$ is the
{\it terminal vertex} map.
\hfill\qed
\end{hypotheses}

We first consider a simple form of retraction, which
amplifies Defini\-tions~III.7.1 of~\cite{DicksDunwoody89}.
Recall that a vertex $v$ of a tree is called a {\it sink} if every 
edge of the tree is oriented towards $v$.

\begin{compression}\label{lem:compress} Suppose that Hypotheses~$\ref{hyp:sec2}$ hold.

Let $E'$ be a $G$-subset of $E$ such that each component of the subforest $T - E'$ of $T$
has a $($unique$)$ sink.  Let $V'$ denote the set of sinks of the components of  $T-E'$.

Let $i \colon E' \to E$ denote the inclusion map,  and 
let $\phi\colon V \to V'$ denote the $G$-retraction which assigns, to each $v \in V$,
the sink of the component of $T - E'$ containing~$v$. 

Then the $G$-graph $T' = (T', V', E', 
\phi \circ \iota \circ i, \phi \circ \tau \circ i)$
is a   $G$-tree.
\end{compression}

Let $E'' = E - E'$ and let $V'' = V - V'$.  
Then $T - E'$ is the $G$-subforest of $T$ with vertex set $V$ and edge set $E''$.
For each $v  \in V$,    $\phi(v)$ is reached in $T$ by starting at $v$
and travelling as far as possible along edges in $E''$ 
respecting the orientation.  The initial vertex map $\iota \colon E \to V$
induces a bijective map $E'' \to V''$.

We say that $T'$ is obtained from $T$ by {\it compressing the 
closures of the elements of
$E''$ to their terminal vertices} or by {\it compressing the components
of $T-E'$ to their sinks}.

In applications, we usually first $G$-equivariantly reorient $T$ and then, 
in the resulting tree,  compress a $G$-set of closed edges 
to their terminal vertices;  we then call the combined  procedure a
 {\it $G$-equivariant compressing operation}.

\begin{proof} [Proof of Lemma~\normalfont\ref{lem:compress}] 
The map $\phi$ induces a surjective $G$-map  $T \to T'$ 
in which the fibres are the components of $T - E'$.  It follows that
$T'$ is a $G$-tree.
\end{proof}

We now recall the sliding operation of Rips-Sela~\cite[p.~59]{RipsSela97}
as generalized by Forester~\cite[Section~3.6]{Forester02}; see 
also the Type 1 operation of~\cite[p.~146]{Dunwoody98}.  
We find it convenient to express the result and the proof
 in the notation of~\cite{DicksDunwoody89}. 

\begin{Forester}\label{lem:sliding}
Suppose that Hypotheses~$\ref{hyp:sec2}$ hold.

Let $e$ and $f$ be elements of  $E$.

Suppose that  $\tau e = \iota f$, $G_e \le G_f$, and  $Gf \cap Ge = \emptyset$.

Let $\tau'\colon E \to V$ denote the map given by 
$$  e' \mapsto \tau'(e'):=
\begin{cases}
\tau(e') &\text{if $e'\in E-Ge$,}\\
\tau(gf) &\text{if $e' = ge$ for some $g \in G$,}
\end{cases}$$
for all $e' \in E$.

Then the $G$-graph $T' = (T', V,E, \iota, \tau')$ is a $G$-tree.
\end{Forester}

Here, we say that 
$T'$ is  obtained from $T$ by {\it $G$-equivariantly sliding  $\tau e$ along $f$ from 
$\iota f$ to $\tau f$}.

In applications, we usually first $G$-equivariantly reorient $Ge$, or $Gf$, or both, or neither,
and then, in the resulting tree,  $G$-equivariantly slide  $\tau e$ 
along $f$ from  $\iota f$ to~$\tau f$, and then reorient
back again.  We then call the combined procedure a {\it $G$-equivariant sliding operation}.

\begin{proof}[Proof of Lemma~$\ref{lem:sliding}$]  It is clear that $T'$ is a $G$-graph.

Let $X$ be the $G$-graph obtained from $T$ by deleting the two edge orbits $Ge \cup Gf$,
and then inserting one new vertex orbit $Gv$ and three new edge orbits $Ge' \cup Gf_1 \cup Gf_2$,
 with $G_{e'} = G_e$, $G_v = G_{f_1} = G_{f_2} = G_f$, and setting $$\iota(e') = \iota(e), \quad
\iota(f_1) = \iota(f) = \tau(e),  \quad \iota(f_2) = \tau(e) = \tau(f_1) = v,  \quad \tau(f_2) = \tau(f).$$ 
Thus we are $G$-equivariantly
 subdividing $f$ into $f_1$ and $f_2$ 
by adding $v$, and then sliding $\tau e$ along $f_1$ 
from $\iota f_1$ to $\tau f_1 = v$.

 Then $T$ is recovered from $X$ by $G$-equivariantly compressing the closure of
$f_1$ to $\iota(f_1)$,  and renaming $f_2$ as $f$, $e'$ as $e$.  Thus 
$X$ maps onto $T$ with fibres which are
trees.  It follows that $X$ is a tree; 
see~\cite[Proposition~III.3.3]{DicksDunwoody89}.  

Also $T'$ is recovered from $X$ by $G$-equivariantly compressing 
the closure of $f_2$ to $\tau(f_2)$,  and renaming $f_1$ as $f$, 
$e'$ as $e$.  By  Lemma~\ref{lem:compress}, $T'$ is a tree.
\end{proof}

\section{Filtrations}

Throughout this section we will be working with the following.

\begin{hypotheses}\label{hyp:1}  Let $T = (T,V,E, \iota,\tau)$ be a  $G$-tree, 
let $U$ be a $G$-retract of the $G$-set~$V$, and let $W = V-U$.
\hfill\qed
\end{hypotheses}

\begin{conventions} 
We shall use interval notation for ordinals; for example,  if $\kappa$ is an ordinal, then
 $[0,\kappa)$  
denotes the set of all ordinals $\alpha$ such that $\alpha < \kappa$.

If we have an ordinal $\kappa$ and a specified map from a set $X$  to $[0,\kappa)$, then
we will understand that the following notation applies.   
Denoting the image of each $x \in X$ by $\degg(x) \in [0,\kappa)$, 
we write,  for each $\alpha \in [0,\kappa)$ and each $\beta \in [0,\kappa]$, 
$$X[\alpha] := \{x \in X \mid \degg(x) = \alpha\} \quad \text{and} \quad
X[0,\beta) := \{x \in X \mid \degg(x) < \beta\}.$$
\vskip -.7cm \hfill\qed \vskip .4cm
\end{conventions}

\begin{definitions}\label{defs:filters} Suppose that Hypotheses~\ref{hyp:1} hold.

Let $P(T)$ denote the set of paths in $T$, 
as in Definitions~I.2.3 of~\cite{DicksDunwoody89}.
Thus, for each $p \in P(T)$, we have the  {\it initial vertex} of $p$, 
denoted $\iota p$, the {\it terminal vertex} of $p$, denoted $\tau p$, 
the {\it set of edges which occur in} $p$,
denoted  $E(p) \subseteq E$, the {\it length} of $p$, denoted $\length(p) \in \naturals$,
and the {\it $G$-stabilizer} of $p$, denoted $G_p \le G$.

Let $\kappa$ be an ordinal and let  
\begin{align}
&T \to [0,\kappa), \quad x \mapsto \degg(x)\label{eq:height1}
 \intertext{be a map. Since $T$ is nonempty, $\kappa$ must be nonzero. As a set, $T = V \cup E$.  Thus,  for each $\alpha \in [0,\kappa)$, 
we have $T[\alpha]$, $E[\alpha]$ and $V[\alpha]$,
and, for each $\beta \in [0,\kappa]$,  we have $T[0,\beta)$, $E[0,\beta)$ and $V[0,\beta)$.
\newline\indent For each $w \in W$, we then define}
&\text{$P_T(w):= \{p \in P(T) \mid
\iota p = w,\, G_{p} = G_{w},\,\degg(\tau p) < \degg(w)$,}\hskip0cm\nonumber \\  
&\hskip4cm \text{$\degg(E(p)) \subseteq \{\degg(w), \, \degg(w)+1\}\}$.}\nonumber
\intertext{\indent We say that~\eqref{eq:height1} is a 
{\it $U$-filtration of $T$} if all of the following hold:}
&\text{for each $\beta \in [0,\kappa]$, $T[0,\beta)$ is a $G$-subforest of~$T$;}\label{1}
\\&\text{$T[0] = U$;}\label{2}
\\&\text{for each $\alpha \in [1,\kappa)$,  $T[\alpha]$  is a $G$-finite $G$-subset of $T$;\, and,}\label{3}
\\&\text{for each $w \in W$, $P_T(w)$ is nonempty.}\label{4}
\end{align}
\vskip -0.8cm \hfill\qed \vskip 0.5cm
\end{definitions}

\begin{lemma}\label{lem:filtration} If Hypotheses~$\ref{hyp:1}$ hold, 
then there exists a $U$-filtration of~$T$.
\end{lemma}

\begin{proof}  We shall recursively construct  
a family $(E[\alpha] \mid \alpha \in [0, \kappa))$
of $G$-subsets of $E$, for some nonzero ordinal $\kappa$.

We take $E[0] = \emptyset$.

Suppose that $\gamma$ is a nonzero ordinal, and that we have 
a family $(E[\alpha] \mid \alpha \in [0, \gamma))$
of $G$-subsets of $E$. 

For each $\beta \in [0,\gamma]$, we define
$$E[0,\beta):= \bigcup\limits_{\alpha \in [0,\beta)} E[\alpha]
\quad\text{and}\quad
V[0,\beta) := \begin{cases}
\emptyset &\text{if $\beta = 0$},\\
U \cup \iota(E[0,\beta)) \cup \tau(E[0,\beta)) &\text{if $\beta > 0 $}.
\end{cases}$$
For each $\alpha \in [0,\gamma)$, we 
 define $V[\alpha] := V[0,{\alpha+1}) - V[0,\alpha)$.
Thus  $$V[0,\beta) = \bigcup\limits_{\alpha \in [0,\beta)} V[\alpha].$$

If $E[0,\gamma) = E$, we take $\kappa = \gamma$ and the construction terminates. 

Now suppose that $E[0,\gamma) \subset E$.  We shall explain how to choose 
$E[\gamma]$.

If $\gamma$ is a limit ordinal or $1$, we take 
 $E[\gamma]$ to be an arbitrary single $G$-orbit in
$E - E[0,\gamma)$. 

If $\gamma$ is a successor ordinal greater than $1$ then there is a unique
$\alpha \in [1,\gamma)$ such that $\gamma = \alpha +1$, and
we want to construct $E[\alpha + 1]$.
Notice that $V[0,\alpha)$ is a $G$-retract of $V$ because $V[0,\alpha)$ contains $U$.
Thus we can $G$-equivariantly specify, for each $w \in V[\alpha]$, 
a $T$-geodesic $p = p(w)$ from $w$ to an element $v = v(w) \in V[0,\alpha)$
fixed by $G_w$.  Since $G_w$ fixes both ends of $p$,
$G_w$ fixes $p$. Hence we may assume that $v$ is the first, and hence only, vertex of $p$ 
that lies in $V[0,\alpha)$. Clearly $G_p$ fixes~$w$.  Thus $G_w = G_p$.
Let $P_{\alpha+1}$ denote the set of edges which occur in the $p(w)$,
as $w$ ranges over $V[\alpha]$.
Then $P_{\alpha+1} \subseteq E - E[0,\alpha)$, 
since each element of $E[0,\alpha)$
has {\it both} vertices in $V[0,\alpha)$.  
If $P_{\alpha+1} \subseteq E[\alpha]$,
we choose $E[\alpha+1]$ to be an arbitrary single $G$-orbit in
$E - E[0,\alpha+1)$.  If $P_{\alpha+1} \not\subseteq E[\alpha]$,
we take $E[\alpha+1] = P_{\alpha+1} - E[\alpha]$.
This completes the description of the recursive construction.

We now verify that we have a $U$-filtration of $T$.

It can be seen that, for each ordinal $\gamma$ such that
$(E[\alpha] \mid \alpha \in [0, \gamma))$
is defined, the $E[\alpha]$, $\alpha \in [1, \gamma)$,
are pairwise disjoint, nonempty, $G$-subsets of $E$.
Hence the cardinal of $\gamma$ is at most one more than the cardinal of $E$.
Therefore the construction terminates at some stage.
This implies that there exists a nonzero ordinal $\kappa$ such that
$E[0,\kappa) = E$. Also $V[0,\kappa) = V$, and 
$(V[\alpha] \mid \alpha \in [0, \kappa))$ 
gives a partition of $V$.  Thus we have an implicit map  
$T \to [0,\kappa)$  and we denote 
it by $x \mapsto \degg(x)$.

Clearly~\eqref{1},~\eqref{2} and~\eqref{4} hold.

If $\alpha \in [1,\kappa)$ and $E[\alpha]$ is $G$-finite, then
either $E[0,\alpha+1) = E$ or
$V[\alpha]$, $P_{\alpha + 1}$ and $E[\alpha +1]$ are 
$G$-finite.  It follows, by transfinite induction, that $E[\alpha]$ 
and $V[\alpha]$ are $G$-finite for all $\alpha \in [1,\kappa)$.
Thus~\eqref{3} holds.
\end{proof}

\section{The main result}
\label{sec:retract}

Let us introduce a technical concept which generalizes 
that of a finite subgroup.

\begin{definitions}
A subgroup $H$ of $G$ is said to be {\it $G$-conjugate incomparable}
if, for each $g \in G$,   $H^g \subseteq H$ (if and) only if $H^g = H$.
This clearly holds if $H$ is finite.

We say that a $G$-set $X$  {\it has $G$-conjugate-incomparable stabilizers}
if, for each $x \in X$, the $G$-stabilizer $G_x$ is a $G$-conjugate-incomparable subgroup,
that is, for each $g \in G$,  $G_x \subseteq G_{gx}$ (if and) only if $G_{x} = G_{gx}$.
\hfill \qed
\end{definitions}

Throughout this section we will be working with the following.

\begin{hypotheses}\label{hyp:2}  Let $T = (T,V,E, \iota,\tau)$ be a  $G$-tree, 
let $U$ be a $G$-retract of the $G$-set~$V$, and let $W = V-U$.

Suppose that the $G$-set  $W$ has $G$-con\-ju\-gate\d1in\-com\-pa\-ra\-ble 
stabilizers.

Let $\kappa$ be an ordinal and let
\begin{equation}\label{eq:height2}
\degg \colon V \cup E \to [0,\kappa), \quad x \mapsto \degg(x), 
\null\hskip 4.5cm
\end{equation}
be a $U$-filtration of $T$. 
 \hfill\qed  
\end{hypotheses}

\begin{definitions}\label{defs:lower}  Suppose that Hypotheses~\ref{hyp:2} hold.

 Let $w \in W$.  Define $d_T(w)  := \min\{\length(p) \mid p \in P_{T}(w)\}$.
Then $d_T(w)$ is a positive integer and
\begin{align}
&\text{$d_T(gw) = d_T(w)$ for all $g \in G$.}\label{eq:Ginvariance}
 \intertext{\indent For $v_0$, $v_1$ in $V$, we say that 
$v_1$ is {\it  lower than} $v_0$ if  
one of the following holds:}
&\text{$\degg(v_0) > \degg(v_1)$;}\label{eq:strong1}\\
&\text{$\degg(v_0) = \degg(v_1) > 0$  and $G_{v_0} < G_{v_1}$; or,}\label{eq:strong2}\\
&\text{$\degg(v_0) = \degg(v_1) > 0$ and  $G_{v_0} = G_{v_1}$  and 
$d_T(v_0) > d_T(v_1)$.}\label{eq:strong3}
\intertext{\indent 
An edge $e$ of $T$ is said to be {\it problematic} 
if it joins vertices $v_0$, $v_1$ such that
$\degg(e) = \degg(v_1) = \degg(v_0)+1$.
Notice that $\degg(e)$ is a successor ordinal
and that $v_0$ is  lower than~$v_1$.
\newline \indent For each $v_0 \in W$, there exists a path}
&\text{$v_0, e_1^{\epsilon_1}, v_1, e_2^{\epsilon_2}, 
v_2, \ldots, e_d^{\epsilon_d}, v_{d}$ in $P_T(v_0)$ such that $d = d_T(v_0)$.}\label{eq:path}
\end{align}
Here $\degg(v_1) \le \degg(v_0) +1$. 
We say that $v_0$ is a {\it problematic} vertex of $T$
if there exists a path as in~\eqref{eq:path} such that $\degg(v_1) = \degg(v_0) +1$. 
In this event $\degg(e_1) = \degg(v_1)$ and $e_1$ is a problematic edge of $T$.
\hfill\qed
\end{definitions}

\begin{lemma}\label{lem:transfinitesliding}
If Hypotheses~$\ref{hyp:2}$ hold, then applying some transfinite sequence of $G$-equivariant 
sliding operations to $T$ yields a $G$-tree $T' = (T',V,E, \iota',\tau')$
such that~\eqref{eq:height2} is also a $U$-filtration of
$T'$ and $T'$ has no problematic vertices.
\end{lemma}

\begin{proof}
We shall construct a family of trees 
$$(T_\beta = (T_\beta, V, E, \iota_\beta, \tau_\beta)\mid \beta \in [0,\kappa])$$
such that, for each $\beta \in [0,\kappa]$,~\eqref{eq:height2} is a $U$-filtration
of  $T_\beta$, and $T_\beta$ has no problematic vertices in $V[0,\beta)$.

We take $T_0 = T$.

For each successor ordinal $\beta = \alpha +1 \in [0,\kappa)$, $T_{\alpha +1}$ will be obtained from
$T_\alpha$ by altering, if necessary, $\iota_\alpha$ and $\tau_{\alpha}$ on 
$E[\alpha+1]$, as described below.

For each limit ordinal $\beta \in [0,\kappa]$, we let
$\iota_\beta$ be given on $E[\alpha]$ by $\iota_\alpha$, for each
$\alpha \in [0,\beta)$,  and similarly for $\tau_\beta$.

Suppose then that $\beta = \alpha +1 \in [0,\kappa)$, that we have a tree
 $T_\alpha = (T_\alpha, V, E, \iota_\alpha, \tau_\alpha)$,
and that \eqref{eq:height2} is a $U$-filtration of $T_\alpha$, 
and that $T_\alpha$ has no problematic vertices in~$V[0,\alpha)$. 

We now describe a crucial {\it problem-reducing procedure} that can be applied
 in the case where there exists some $v_0 \in V[\alpha]$
which is a problematic vertex of  $T_\alpha$.  

Let $d = d_{T_\alpha}(v_0)$. 
Thus, there exists a path
$$v_0, e_1^{\epsilon_1}, v_1, e_2^{\epsilon_2}, 
v_2, \ldots, e_d^{\epsilon_d}, v_{d}$$
in $P_{T_\alpha}(v_0)$ such that
$v_1 \in V[\alpha+1]$.  Hence,  $e_1 \in E[\alpha+1]$.
Without loss of generality, let us assume that $\epsilon_1 = -1$.

There exists a least $i \in [2,d]$ 
 such that $v_i \in V[0,\alpha +1)$.  
Then $$\{v_1,\ldots, v_{i-1}\} \subseteq V[\alpha+1] \quad \text{and, hence,} \quad 
\{e_1,\ldots, e_{i}\} \subseteq E[\alpha+1].$$

We claim that 
 $Ge_1 \,\cap\, \bigcup\limits_{j=2}^i Ge_j = \emptyset$.  
Suppose this fails.  Then $e_1 \in \bigcup\limits_{j=2}^i Ge_j$.  Here,
  $v_0 \in \bigcup\limits_{j=1}^i Gv_j$.
Since $v_0 \in V[\alpha]$ and  
$\bigcup\limits_{j=1}^{i-1} Gv_j \subseteq V[\alpha+1]$
we see that $v_0 \in Gv_i$.  Hence $v_i \in V[\alpha]$ and, 
by~\eqref{eq:Ginvariance},  $d_{T_\alpha}(v_i) = d_{T_\alpha}(v_0) = d$.
But $G_{v_0} = G_p \subseteq G_{v_i}$. 
Since $G_{v_0}$ is a $G$-con\-ju\-gate-incomparable subgroup, $G_{v_0} = G_{v_i}$.
It follows that 
$$v_i, e_{i+1}^{\epsilon_{i+1}}, v_{i+1}, \ldots, e_d^{\epsilon_d}, v_{d}$$
lies in $P_{T_\alpha}(v_i)$.  Hence $d_{T_\alpha}(v_i) \le d-i$, which is a contradiction.
This proves the claim.   

 By Lemma~\ref{lem:sliding}, we
can $G$-equivariantly slide $\iota e_1$ along $e_2^{\epsilon_2}$ 
from $v_1$ to $v_2$, and then
$G$-equivariantly slide $\iota e_1$ along $e_3^{\epsilon_3}$ 
from $v_2$ to $v_3$, and so on, up to $v_i$.  We then get a new $G$-tree 
$T_{\alpha,1} = (T_{\alpha,1},V,E,\iota_{\alpha,1},\tau_{\alpha,1})$ by
 $G$-equivariantly sliding  $\iota e_1$ along our path from $v_1$ to $v_i$.

Let $e_1'$ denote $e_1$ viewed as an edge of $T_{\alpha,1}$. Wherever $v_1,e_1,v_0$ 
occurs in a path in $T_\alpha$, it can be replaced with the sequence  
$$v_1,e_2^{\epsilon_2},v_2, \ldots,  v_{i-1}, 
e_i^{\epsilon_i}, v_i, e_1^{\prime}, v_0$$ to obtain a path in $T_{\alpha,1}$.  
It is important to note that all the edges
involved here lie\linebreak in~$E[\alpha +1]$.
In terms of  the free groupoid on $E[\alpha +1]$,
$e_1 = e_2^{\epsilon_2} e_3^{\epsilon_3}\cdots e_i^{\epsilon_i} e_1'$, and
 we are performing the
change-of-basis which replaces $e_1$ with $e_1'$.

It is easy to see that~\eqref{1}--\eqref{4}
then hold for $T_{\alpha,1}$.  Thus~\eqref{eq:height2} is a $U$-filtration of $T_{\alpha,1}$. 
 Notice that $T_{\alpha,1}$, like $T_\alpha$, has no problematic vertices in $V[0,\alpha)$.
We have reduced the number of $G$-orbits of problematic edges in $E[\alpha+1]$.  

This completes the description of a problem-reducing procedure.

Since  $E[\alpha+1]$ is $G$-finite by~\eqref{3}, on repeating problem-reducing procedures as often as possible, 
we find some $m \in \naturals$, and a sequence
$$T_\alpha = T_{\alpha,0},\, T_{\alpha,1},\, \ldots,\, T_{\alpha,m},$$
such that $T_{\alpha,m}$ has no problematic vertices in $V[0,\alpha) \cup V[\alpha] = V[0,\alpha+1)$.
We define $T_{\alpha+1} = (T_{\alpha+1}, V,E,\iota_{\alpha+1}, \tau_{\alpha+1})$ to be 
$T_{\alpha,m}$.  Notice that $\iota_{\alpha+1}$ agrees with 
$\iota_\alpha$ on $E - E[\alpha+1]$, and similarly for $\tau_{\alpha+1}$.

Continuing this procedure transfinitely, we arrive at a tree $T_\kappa$ which has no 
problematic vertices.
\end{proof}

\begin{lemma}\label{lem:compresstoU} If Hypotheses~$\ref{hyp:2}$ hold and $T$ has no problematic vertices, 
then applying some $G$-equivariant compressing operation on $T$ yields a $G$-tree 
with vertex set~$U$.
\end{lemma}

\begin{proof}  We claim that any sequence in $V$ is finite if each term 
is  lower than all its predecessors.

Let $\alpha \in [0,\kappa)$.

If $v_0$, $v_1$ are elements of the same $G$-orbit of $V[\alpha]$, 
then $v_1$ is  not lower than $v_0$, that is,~\eqref{eq:strong1}--\eqref{eq:strong3} all fail;
this follows from~\eqref{eq:Ginvariance} and the fact that $V[\alpha]$ has
$G$-conjugate-incomparable stabilizers.  

Thus, if $n \in \naturals$ and $v_1, v_2, \ldots, v_{n}$ 
is a sequence in $V[\alpha]$ such that each term 
is  lower than all its predecessors, then $Gv_1, Gv_2, \ldots, Gv_n$ 
are pairwise disjoint,  and $n$ is at most the number of $G$-orbits in $V[\alpha]$. 
It follows that any sequence in $V[\alpha]$ is finite if each term 
is  lower than all its predecessors.  The claim now follows.

Let us $G$-equivariantly reorient $T$ so that,
for each edge $e$, $\iota e$ is not lower than~$\tau e$. 

Let $v_0 \in W$.  Let us $G$-equivariantly choose  a path  
$$v_0, e_1^{\epsilon_1}, v_1, e_2^{\epsilon_2}, 
v_2, \ldots, e_d^{\epsilon_d}, v_{d}$$
in $P_{T}(v_0)$ such that  $d = d_T(v_0)$.  Then  
we call $e_1$ the {\it distinguished edge} associated to~$v_0$, and 
$v_1$ the {\it distinguished neighbour} of~$v_0$. 

Let $E''$ denote the set of distinguished edges chosen in this way.

Let us consider the above path for $v_0$. From Definitions~\ref{defs:lower},  
we see that, since $T$ has no problematic vertices, $\degg(v_0) \ge \degg(v_1)$.
We claim that  $v_1$ is  lower than~$v_0$.  The claim is clear if
 $\degg(v_0) > \degg(v_1)$ (in which case, $d=1$), and we may assume that
$\degg(v_0) = \degg(v_1)\,\, (>0)$.  
Again, the claim is clear if $G_{v_0} < G_{v_1}$, and
we may assume that $G_{v_0} = G_{v_1}$.  Here $G_{v_1}$ fixes $p$, and
the path  $$v_1, e_2^{\epsilon_2}, v_2, \ldots, e_d^{\epsilon_d}, v_{d}$$
shows that $d_T(v_1) \le d-1 < d =  d_T(v_0)$, and the claim is proved.
Hence $\epsilon_1 = 1$.  

Thus  $\iota$ induces a bijection $E'' \to W$.

Moreover, in travelling along the distinguished edge $e_1$
respecting the orientation, from  $v_0$ to its distinguished neighbour $v_1$, we move
 to a  lower vertex.

Thus, starting at any element $v$ of $V$, 
after travelling a finite number of steps along distinguished 
edges respecting the orientation, we arrive
at a vertex, denoted $\phi(v)$,  
with no distinguished neighbours, that is,   $\phi(v) \in U$.
 
By Lemma~\ref{lem:compress}, compressing the closures of the 
distinguished edges to their terminal vertices
gives a $G$-tree with vertex set $U$ 
and edge set $E - E''$.
\end{proof}

We now come to our main result.  In Section~\ref{sec:example}, 
we will see that
the $G$-con\-ju\-gate\d1in\-com\-parability hypotheses cannot be omitted.

\begin{theorem}\label{thm:main2}
Let  $T$ be a $G$-tree,
and let $U$ be a $G$-retract of the $G$-set~$VT$.
Suppose that the $G$-set $ET$ has $G$-con\-ju\-gate\d1in\-com\-pa\-ra\-ble 
stabilizers, or, more generally, that the  $G$-set  
$VT-U$ has $G$-con\-ju\-gate\d1in\-com\-pa\-ra\-ble 
stabilizers.

Then applying to $T$ some transfinite sequence of $G$-equivariant sliding
operations followed by some $G$-equivariant compressing operation
yields a $G$-tree $T'$  such that\linebreak $VT' = U$.

Here $ET'$ is a $G$-subset of $ET$, and there exists a $G$-set isomorphism $$ET - ET' \simeq VT - VT' = VT - U.$$
\end{theorem}

\begin{proof}  For each $w \in VT-U$, there exists $u \in U$
such that $G_w \le G_u$.  If $e$ denotes the first edge in the $T$-geodesic
from $w$ to $u$, then $G_e = G_w$.  Thus, if $E$ has $G$-conjugate-incomparable 
stabilizers, then the same holds for $VT-U$.

By Lemma~\ref{lem:filtration}, we may assume that Hypotheses~\ref{hyp:2} hold.
By Lemma~\ref{lem:transfinitesliding}, 
we may assume that
$T$  itself has no problematic vertices.
Applying Lemma~\ref{lem:compresstoU}, we obtain the result;
the final assertion follows from Lemma~\ref{lem:compress}.
\end{proof}

We record the special case of
Theorem~\ref{thm:main2} that is of interest to us.

\begin{retraction}\label{lem:main}
Let $T$ be a $G$-tree whose edge stabilizers are finite,
and let $U$ be any $G$-retract of the $G$-set $VT$. 
Then there exists a $G$-tree whose edge stabilizers are finite 
 and whose vertex set is the $G$-set $U$. 
\hfill\qed
\end{retraction}

\section{The almost stability theorem and applications}
\label{sec:applications}

We now combine the almost stability theorem~\ref{thm:almost} and the retraction lemma~\ref{lem:main}. 

\begin{theorem}\label{thm:main}  Let $E$ and $A$ be $G$-sets such that
 $E$ has finite stabilizers and $A$ has trivial $G$-action.
If  $V$ is a $G$-retract of a $G$-stable almost equality class in 
$(E,A)$,  then there exists a $G$-tree whose edge stabilizers are finite
 and whose vertex set is the $G$-set~$V$.
\end{theorem}

\begin{proof} By the almost stability theorem~\ref{thm:almost}, there exists a $G$-tree
whose edge stabilizers are finite and whose vertex set is the given 
$G$-stable almost equality class in $(E,A)$.  
By the retraction lemma~\ref{lem:main},
there exists a $G$-tree whose edge stabilizers are
 finite and whose vertex set is $V$.  
\end{proof}

We now recall Definitions~IV.2.1 and IV.2.2 of~\cite{DicksDunwoody89}.

\begin{definitions} Let $M$ be a {\it $G$-module}, that is, an additive abelian 
group which is also a $G$-set such that $G$ acts as group automorphisms on $M$.
Thus a $G$-module is simply a left module over the integral group ring $\integers G$.

If $d \colon G \to M$ is a {\it derivation}, that is, a map such that
$d(xy) = d(x) + xd(y)$ for all $x$, $y \in G$, then
$M_d$ denotes the set $M$ endowed with the $G$-action 
$$G \times M \to M, \quad (g,m) \mapsto  g \cdot m := gm + d(g)
\quad \text{for all $g\in G$ and all $m \in M$.}$$ 
It is straightforward to show that $M_d$ is a $G$-set.
This construction has  
made other appearances in the 
literature; see~\cite[Remarque~4.a.5]{delaHarpeValette89}.

We say that  $M$ is  an {\it induced} 
$G$-module if there exists an abelian group $A$ such that 
$M$ is isomorphic, as $G$-module, to $AG := \integers G \otimes _\integers A$.

We say that $M$ is a {\it $G$-projective} $G$-module if 
$M$ is isomorphic, as $G$-module, to a direct summand of an induced $G$-module.
\hfill\qed
\end{definitions}

\begin{example} If $R$ is any ring and $P$ is a projective left $RG$-module, then
there exists a free left $R$-module $F$ such that $P$ is isomorphic, as $RG$-module, 
to an $RG$-summand of 
$$RG \otimes _R F = \integers G \otimes _\integers R \otimes_R F = \integers G \otimes _\integers F =  FG .$$
Hence $P$ is $G$-projective.
\hfill\qed
\end{example}

The following generalizes Theorem~IV.2.5 and Corollary~IV.2.8 of~\cite{DicksDunwoody89}.

\begin{theorem}\label{thm:derivation} If $P$ is a $G$-projective $G$-module, and $d\colon G \to P$ is
a derivation, then there exists a $G$-tree
 whose edge stabilizers are finite and whose vertex set is the $G$-set $P_d$.
\end{theorem}

\begin{proof} There exists an abelian group $A$ such that $P$ is isomorphic to 
a $G$-summand of $AG$.
We view $P$ as a $G$-submodule of $AG$.  There exists an additive 
$G$-retraction $\pi \colon AG \to P$.

We view $AG$ as the almost equality class of $(G,A)$ which contains the zero map.  
Thus $AG$ is a $G$-submodule of $(G,A)$, and we have a derivation
 $$d:G \to P \subseteq AG \subseteq (G,A).$$ 
 
By a classic result of Hochschild's,  
there exists $v \in (G,A)$ such that, for all $g \in G$, $d(g) = gv -v$.
For example,  we can take $v\colon x \mapsto -(d(x))(x)$, for all $x \in G$. 
See the proof of Proposition~IV.2.3 in~\cite{DicksDunwoody89}.

Let $U = v +P$ and $V = v+AG$.  Then $U \subseteq V \subseteq (G,A)$, 
and $V$ is the almost equality class which contains $v$.  Also, $U$ and $V$ are
$G$-stable, since, for each $g \in G$, $gv = v + d(g) \in v + P \subseteq v + AG$.  
The map
 $$V \to U, \quad v+m \mapsto v +\pi(m), \text{ for all } m \in AG,$$
is a $G$-retraction, since, for all $m \in AG$,
\begin{align*}
g(v+m) = v + gm + d(g) \quad \mapsto \quad v + \pi(gm + d(g)) 
&=  v + g\pi(m) + d(g) \\&= g(v+\pi(m)).
\end{align*}

By Theorem~\ref{thm:main}, 
there exists a $G$-tree
whose edge stabilizers are finite and whose vertex set is the $G$-set $U$. 
 
The bijective map $P \to U$, $p \mapsto v + p$, is an isomorphism
of $G$-sets $P_d \xrightarrow{\null_\sim} U$.    Now the result follows.
\end{proof}

\begin{remark} Notice that, in Theorem~\ref{thm:derivation}, the stabilizer
of a vertex $p\in P_d$ is precisely the kernel of the derivation 
$$d + \ad p\colon  G \to P, \quad g \mapsto d(g) + gp -p = (g-1)(v + p).$$
\null\vskip -1.2cm \hfill\qed \vskip 0.1cm
\end{remark}

The following generalizes Corollary~IV.2.10 of~\cite{DicksDunwoody89}
and is used in the proof of Lemma~5.16 of~\cite{DicksLinnell06}.

\begin{corollary}\label{cor:coset} Let $M$ be a $G$-module, let $P$ be a 
$G$-projective $G$-submodule of~$M$,  and let $v$ be an element of $M$.  
If the subset $v + P$ of $M$ is $G$-stable, 
then there exists a $G$-tree whose edge stabilizers are finite and 
whose vertex set is the $G$-set  $v +P$.
\end{corollary}

\begin{proof} The inner derivation $\ad  v \colon G \to M$ restricts to 
a derivation $d \colon G \to P$, $g \mapsto gv - v \in P \subseteq M$,
for all $g \in G$.  The  bijective map $P  \to v + P$, $p \mapsto v+p$, is 
then an isomorphism of $G$-sets $P_d \xrightarrow{\null_{\sim}} v+P$.  
Now the result follows from 
Theorem~\ref{thm:derivation}.
\end{proof}

\begin{example}  Let $R$ be a nonzero associative ring, and let $\omega R G$ be the
augmentation ideal of the group ring $RG$.  

Notice that, in the (left) $G$-set $RG$,  both the coset $1 + \omega R G$ and
$RG -\{0\}$ are $G$-stable, and that the  $G$-set 
$RG -\{0\}$ has finite stabilizers.

If $\omega R G$ is projective as left $RG$-module, then, by
Corollary~\ref{cor:coset}, there exists  a $G$-tree $T$ with 
$VT = 1 + \omega R G \subseteq RG - \{0\}$; hence $T$ has finite stabilizers. 
This sheds some light on the main step in the characterization of 
 groups of cohomological dimension at most one over $R$.
See, for example,~\cite[Theorem IV.3.13]{DicksDunwoody89}.
\hfill\qed
\end{example}

\section{A more general form}
\label{sec:generalized}

We next want to generalize Theorem~\ref{thm:main}. 

The following is similar to Lemma~2.2 of~\cite{DicksKropholler95},
and the proof is straightforward.

\begin{lemma}\label{lem:iso} Let $E$ and $A$ be $G$-sets such that, for each $e \in E$, 
$G_e$  acts trivially on~$A$.

Let $\bar A$ denote the $G$-set with the same underlying set as
$A$ but with trivial $G$-action.     

Let $E_0$ be a $G$-transversal in $E$.  

For each $\phi \in (E,A)$, 
let $\widehat \phi \in (E,\bar A)$ be defined by $\widehat \phi(ge) = g^{-1} \cdot \phi(ge)$
for all $(g,e) \in G \times E_0$, where $\cdot$ denotes the $G$ action on $A$.  

For each $\psi \in (E,\bar A)$, 
let $\widetilde \psi \in (E,A)$ be defined by $\widetilde \psi(ge) =  g \cdot \psi(ge)$
for all $(g,e) \in G \times E_0$.  

Then  
$$(E,A) \to (E,\bar A), \quad \phi \mapsto \widehat \phi, 
\quad  \text{and} \quad (E,\bar  A) \to (E,A),\quad \psi \mapsto \widetilde \psi,$$
are mutually inverse isomorphisms of $G$-sets which preserve 
almost equality between functions.
\hfill\qed
\end{lemma}

Combined, Lemma~\ref{lem:iso} and Theorem~\ref{thm:main}  give the 
most general form  that we know of the almost stability theorem.

\begin{theorem}\label{thm:free}  Let $E$ and $A$ be $G$-sets such that, for each $e \in E$, 
$G_e$ is finite and acts trivially on $A$.
If $V$ is a $G$-retract of a $G$-stable almost equality class
in $(E,A)$, 
then there exists a $G$-tree whose edge stabilizers are finite and
whose vertex set is the $G$-set $V$.
\hfill\qed
\end{theorem}

 For each $e \in E$, if $G_e$ is trivial, then $G_e$ is finite and acts trivially on $A$.
It was this case that was useful in~\cite{DicksKropholler95}.

\section{An example}
\label{sec:example}

In this section, we shall give an example of a group $G$ and a retract of a 
vertex set of a $G$-tree that is not the vertex set of any $G$-tree.

We shall use two technical lemmas.
Recall that, for $x, \,y\in G$, $x^y$ denotes $y^{-1}xy$.

\begin{lemma}\label{lem:schreier}  Let $G = \gen{x,y \mid \quad}$, 
let $n \in \naturals$, 
 and let $g \in G$.
\begin{enumerate}[\normalfont(i)]
\vskip-0.6cm \null 
\item\label{Ione} If  
$x^{2^n}y^{2^n}x^{2^n} \in \gen{x^{2},  y^{2}}^g$,  then\,\,
$n \ne 0$ and $g \in \gen{x^{2},  y^{2}}$.
\vskip-0.6cm \null
\item\label{Itwo} If  
$x^{2^n}y^{2^n}x^{2^n} \in \gen{x^4, xyx, y^4}^g$,  then\,\,
$n \ne 1$ and $g \in \gen{x^4, xyx, y^4}$.
\end{enumerate}
\end{lemma}

\begin{proof}  Let $T = X(G,\{x,y\})$, the Cayley graph of $G$ with respect to $\{x,y\}$, as\linebreak 
in~\cite[Definitions~I.2.1]{DicksDunwoody89}.  Each (oriented) edge of $T$ 
is labelled $x$ or $y$.

Let $H \le G$, and let $w =x^{2^n}y^{2^n}x^{2^n} \in G$.
  Let $X := H\backslash T$, let $Y:= \gen{w}\backslash T$, and let $Z:= G \backslash T$.  

The pullback of the two natural maps $X \to Z$, $Y \to Z$ provides 
detailed information about all nontrivial subgroups of $G$
of the form $\gen{w} \cap H^g$; see~\cite[p.~380]{Dicks94}.  However, this pullback can be
rather cumbersome and we do not require detailed
information.  For our purposes, special considerations will 
suffice, as follows.   

Define $g^{-1}X := (H^g)\backslash T$. 

There is a graph isomorphism $X \simeq g^{-1}X$, $Hx \leftrightarrow H^g g^{-1}x$.

The fundamental group of $X$ with basepoint $H1$, $\pi(X, H1)$, is
naturally isomorphic to $H$, with the elements of $H$ being read off 
closed paths based at~$H1$.  

Similarly,  $H^g$ is naturally isomorphic to $\pi(g^{-1}X, H^g1)$, 
and this in turn is naturally isomorphic to  $\pi(X, Hg)$ via the
graph isomorphism $g^{-1}X \simeq X$.

Suppose that $w$ lies in $H^g$. 
Then $w$ can be read off a closed path in $X$ based at~$Hg$.  
Since $w$ is a cyclically reduced word, the closed path is
cyclically reduced.  
The smallest subgraph of $X$ which contains 
all the cyclically reduced closed paths in $X$ is
called the {\it core} of~$X$, denoted $\core(X)$.  
It follows that the vertex $Hg$ lies in $\core(X)$, and 
that we can start at $Hg$, read $w$ and stay inside
$\core(X)$.

(i) Suppose that $H = \gen{x^2,y^2}$. 

Here $\core(X)$ has  vertex set 
 $\{H1, \, Hx,\, Hy\}$ and labelled-edge set 
$$\{(H1,x,Hx), (Hx,x, Hx^2), (H1, y, Hy), (Hy, y, Hy^2)\}$$
with $Hx^2 = Hy^2 =H1$.  

We note that $Hxy$ and $Hyx$ are outside $\core(X)$.

Since $(Hy)x = Hyx$ does not lie in $\core(X)$, we see that
$Hg \ne Hy$.  Hence,  $Hg \in \{H1,Hx\}$.

Notice that $(H1)(xy) = Hxy$  and $(Hx)(xyx) = Hyx$.
These lie outside $\core(X)$. 
Thus $n \ne 0$. Hence, $x^{2^n} \in H$.

Notice that $(Hx)(x^{2^n}y) = Hxy$  lies outside  $\core(X)$.
Thus $Hg \ne Hx$.  Hence, $Hg = H1$, that is, $g \in H$.

This proves (i).

(ii). Suppose that $H = \gen{x^4, xyx, y^4}$.

Here   $\core(X)$ has   vertex set  
$$\{H1\} \cup \{Hx^i,\, Hy^i \mid 1 \le i \le 3\}.$$ 
and labelled-edge set
$$\{(Hx^{i},x,Hx^{i+1}),\, (Hy^{i},y,Hy^{i+1}) \mid 0 \le i \le 3\} \cup \{(Hx,y, Hxy)\},$$
with $Hx^4 = Hy^4 = H1$  and $Hxy = Hx^3$. 

We note that $Hxy^2 = Hx^3y$, $Hx^2y$,  $Hyx$, $Hy^2x$ and $Hy^3x$,  
all lie outside $\core(X)$.

For any $j$ with $1 \le j \le 3$, 
$(Hy^j)(x) = Hy^jx$ lies outside  
$\core(X)$.  It follows that  $Hg \ne Hy^j$. Hence $Hg = Hx^i$ for some $i$ with $0 \le i \le 3$.

Notice that $(Hx)(xy) = Hx^2y$, $(Hx^2)(xy) = Hx^3y$,
and $(Hx^3)(xyx) = Hyx$.  These all lie outside $\core(X)$. 
Thus, if $n = 0$, then $Hg = H1$. 

Notice that $(H1)(x^2y) = Hx^2y$, $(Hx)(x^2y) = Hx^3y$,
$(Hx^2)(x^2y^2x) = Hy^2x$, and  $(Hx^3)(x^2y^2) = Hxy^2$.
These all lie outside $\core(X)$.  Thus $n \ne 1$.

Now suppose that $n \ge 2$.  Thus $x^{2^n} = (x^4)^{2^{n-2}} \in H$.

Notice that $(Hx)(x^{2^n}y^2) = Hxy^2$, 
$(Hx^2)(x^{2^n}y) = Hx^2y$, 
and  $(Hx^3)(x^{2^n}y) = Hx^3y$.
These all lie outside $\core(X)$.   Thus $Hg = H1$.

This proves (ii).
\end{proof}

It is straightforward to prove the following.

\begin{lemma}\label{lem:really} Let
$G =  \gen{x ,y, t \mid  x^{4t} = x^8,\,  y^{4t} = y^8,\, x^{t^2}y^{t^2}x^{t^2} = 
x^4y^4x^4}$  and let $n \in \naturals$.
\begin{enumerate}[\normalfont (i)]
\vskip-0.5cm \null
\item\label{it:nnotone}
If $n \ne 1$, then $(xyx)^{t^n} = x^{2^n}y^{2^n}x^{2^n}$ in $G$.
\vskip-0.5cm \null
\item\label{it:narb}
$(xyx)^{t^{n+2}} 
= (x^4)^{2^n}(y^4)^{2^n}(x^4)^{2^n}$ in $G$.\hfill\qed 
\end{enumerate}
\end{lemma}

Throughout the remainder of the section we work with the following example.

\begin{hypotheses}\label{hyp:example}  Let
$G =  \gen{x ,y, t \mid  x^{4t} = x^8,\,  y^{4t} = y^8,\, x^{t^2}y^{t^2}x^{t^2} = 
x^4y^4x^4}.$

Let $T = (T,V,E,\iota, \tau)$ be the $G$-graph given by the following data, where $\vee$ denotes
the disjoint union:
$$V = Gu \vee Gw, \quad G_{u} = \gen{x,\,y}, \quad G_{w} =  \gen{x^{4}, \, y^{4}},$$
$$E = Ge \vee Gf, \quad G_{e} = \gen{ x^4,\, xyx,\, y^4}, 
\quad G_{f} = \gen{x^{4}, \, y^{4}},$$
$$\iota(e) = u, \quad \tau(e) = t^2w, \quad
\iota(f) = w, \quad \tau(f) = tw.$$
Using Lemma~\ref{lem:really}, we see that the following hold:
$$G_{e} \le G_{u},\quad  
G_{t^{-2}e} = G_{e}^{t^2} = \gen{x^{16},\,x^4y^4x^4,\, y^{16}} \le G_{w},$$
$$G_{f}= G_{w}, \quad 
G_{t^{-1}f} = G_{f}^{t} = \gen{x^{8},\,y^{8}} \le G_{w}.$$
Thus $T$ is a well-defined $G$-graph.

Let $U = Gu$.  

Let $H = \gen{x,y} \le G$.

For any subset $S$ of $T$, we let $S^{xyx}$ denote $\{s \in S \mid (xyx)s = s\}$.
\hfill\qed
\end{hypotheses}

Since $G_w \le G_u$, it is clear that $U$ is a $G$-retract of $V$.
We shall see that $T$ is a $G$-tree, and that no $G$-tree has vertex set $U$.

\begin{lemma}\label{lem:bass-serre} If Hypotheses~$\ref{hyp:example}$ hold,
then the $G$-graph $T$ is a tree, and $H$ is freely generated by $\{x,y\}$.
\end{lemma}

\begin{proof} Let us momentarily forget Hypotheses~$\ref{hyp:example}$.

Let $Y= (Y, \overline V, \overline E, \overline \iota, \overline \tau)$ 
be the graph given as follows.
$$\overline V = \{\overline u,\overline w\}, \quad 
\overline E = \{\overline e, \overline f\}, \quad 
\overline \iota(\overline e) = \overline u, \quad 
\overline \tau(\overline e) = \overline \iota(\overline f) 
=  \overline \tau(\overline f) = \overline w.$$
Let $Y_0 := (Y_0, \overline V, \{\overline e\}, \overline \iota, \overline \tau)$ 
be the unique maximal subtree of $Y$.

Using the notation of Definitions~I.3.1 of~\cite{DicksDunwoody89},
let $(G(-),\,Y)$ 
be the graph of groups given by the following data.
$$G(\overline u) = \gen{x,\, y \mid \quad}, \quad 
G(\overline w) = \gen{x',\,y' \mid \quad}, \quad
G(\overline e) = \gen{x^4,\, xyx,\, y^4},  \quad 
G(\overline f) = \gen{x',\,y'},$$ 
$$(x^4)^{t_{\overline e}} = x^{\prime 4}, \quad (x yx)^{t_{\overline e}} = x'y'x', \quad  
(y^4)^{t_{\overline e}}  = y^{\prime 4}, \quad (x')^{t_{\overline f}} = x^{\prime 2}, 
\quad (y')^{t_{\overline f}} = y^{\prime 2}.$$ 
Recall that, in the notation of Definitions~I.3.1 of~\cite{DicksDunwoody89}, $(-)^{t_{\overline e}}$ denotes the
edge-group monomorphism associated to $\overline e$.

Let $G := \pi(G(-),Y,Y_0)$, as in Definitions~I.3.4 of~\cite{DicksDunwoody89}.  
Writing $t$ for the element of $G$ that realizes the monomorphism 
$t_{\overline f}\colon G(\overline f) \to G(\overline w)$, 
we have  $$G =   \gen{x, y, x', y', t \mid x^4 = x^{\prime 4}, \, xyx = x'y'x', \, y^4 = y^{\prime 4}, \,
x^{\prime t} = x^{\prime 2},\, y^{\prime t} = y^{\prime 2}}.$$

Then $\gen{x,y \mid \quad} = G(\overline u) \le G$ by
Corollary~I.7.5 of~\cite{DicksDunwoody89}.

Now $x'^{t^2} = x'^{2t} = x^{\prime 4} = x^4$. Thus
$x' = x^{4t^{-2}}$.  Similarly,  $y' = y^{4t^{-2}}$. 
Hence we can write
\begin{align*}
G &= \gen{x,y,t\mid  x^4 = x^{16t^{-2}},\, xyx = x^{4t^{-2}}y^{4t^{-2}}x^{4t^{-2}},\,
 y^4 = y^{16t^{-2}},\\
&\hskip 1.8cm  x^{4t^{-1}}= x^{8t^{-2}},\hskip 3.4cm  y^{4t^{-1}} = y^{8t^{-2}} }\\
&= \gen{x,y,t\mid   x^{4t^2}= x^{16},\,  x^{t^2}y^{t^2}x^{t^2}
= x^{4}y^{4}x^{4},\, y^{4t^2} = y^{16},\,\\
&\hskip 1.8cm x^{4t} = x^{8},\hskip 3.4cm  y^{4t} = y^{8} }\\
&= \gen{x,y,t\mid   x^{4t} = x^{8},\,    x^{t^2}y^{t^2}x^{t^2}= x^{4}y^{4}x^{4},\,
 y^{4t} = y^{8} }.
\end{align*}

Let $T = (T,V,E,\iota,\tau)$ be $T(G(-),Y, Y_0)$, as in Definitions~I.3.4 of~\cite{DicksDunwoody89}.
Thus  
$$V = G \overline u \vee G \overline w, \quad 
G_{\overline u} = \gen{x, \, y}, \quad 
G_{\overline w} = \gen{x',\,y'} = \gen{x^4,\,y^4}^{t^{-2}},$$
$$E = G \overline e \vee G \overline f, \quad 
G_{\overline e} = \gen{x^4,\, xyx,\, y^4},  \quad 
G_{\overline f} = \gen{x',\,y'} = \gen{x^4,\,y^4}^{t^{-2}},$$
$$\iota(\overline e) = \overline u, \quad 
\tau(\overline e) = \overline w, \quad
\iota(\overline f) = \overline w, \quad 
\tau(\overline f) = t \overline w.$$
By Bass-Serre Theory, $T$ is a $G$-tree; see~\cite[Theorem I.7.6]{DicksDunwoody89}.

Let $u:= \overline u$,  $w := t^{-2}\overline w$, 
$e:= \overline e$, $f:= t^{-2}\overline f$.

Then $\iota e = u$, $\tau e = t^2w$, $\iota f = w$, $\tau f = tw$.

Thus the above $G$ and $T$ agree with the $G$ and $T$ of Hypotheses~\ref{hyp:example}, and the result is proved.
\end{proof}

\begin{lemma}\label{lem:prefixed} Let $n \in \naturals$. 
If Hypotheses~$\ref{hyp:example}$ hold, then the following also hold. 
\begin{enumerate} [\normalfont (i)]
\vskip-0.4cm \null
  \item  $(t^nG_ue)^{xyx} = \{t^ne\}$ if $n \ne 1$.
\vskip-0.4cm \null
  \item  $(t^{n+2}G_wt^{-2}e)^{xyx} = \begin{cases}
\{t^{n}e\} &\text{if $n \ne 1$,}\\
\emptyset &\text{if $n = 1$}.
\end{cases}$
  \item  $(t^{n+2}G_wt^{-1}f)^{xyx}  = \begin{cases}
\{t^{n+1}f\} &\text{if $n \ne 0$,}\\
\emptyset &\text{if $n = 0$}.\end{cases}$
  \item  $(t^{n+2}G_wf)^{xyx} = \{t^{n+2}f\}$.
\end{enumerate}
\end{lemma}

\begin{proof} (i). Let $g \in G_u =  \gen{x,y}$.

 Suppose that $n \ne 1$ and that $(xyx)t^nge = t^nge$.  Then $(xyx)^{t^ng} \in G_e.$
By Lemma~\ref{lem:really}\eqref{it:nnotone}, 
$$(x^{2^n}y^{2^n}x^{2^n})^g \in G_e = \gen{x^4,\, xyx,\, y^4}.$$
By Lemma~\ref{lem:schreier}\eqref{Itwo}, 
$g \in \gen{x^4, xyx, y^4}  = G_e$.
Hence $t^nge = t^ne$. 
It is now easy to see that (i) holds.

(ii). Let $g \in G_w =  \gen{x^4,y^4}$.

Suppose that $(xyx)t^{n+2}gt^{-2}e\!=\!t^{n+2}gt^{-2}e$.  Then
$(xyx)^{t^{n+2}gt^{-2}} \in G_e$. By Lemma~\ref{lem:really}\eqref{it:narb},  
$$((x^4)^{2^n}(y^4)^{2^n}(x^4)^{2^n})^g \in G_e^{t^2} = \gen{x^4, xyx, y^4}^{t^2} 
= \gen{x^{16}, x^4y^4x^4, y^{16}}.$$
By Lemma~\ref{lem:schreier}\eqref{Itwo},  $n \ne 1$ and 
$g \in \gen{x^{16}, x^4y^4x^4, y^{16}} = G_e^{t^{2}}$.
Hence $t^{n+2}gt^{-2}e = t^{n}e$.  It is now clear that (ii) holds.

(iii). Let $g \in G_w =  \gen{x^4,y^4}$.

Suppose that $(xyx)t^{n+2}gt^{-1}f\!=\!t^{n+2}gt^{-1}f$.
Then $(xyx)^{t^{n+2}gt^{-1}}\!\in\!G_f$. By Lemma~\ref{lem:really}\eqref{it:narb},  
$$((x^4)^{2^n}(y^4)^{2^n}(x^4)^{2^n})^g \in G_f^{t} = \gen{x^4, y^4}^t = \gen{x^8, y^8}.$$
By Lemma~\ref{lem:schreier}\eqref{Ione},  
$n \ne 0$ and $g \in \gen{x^8, y^8} = G_f^{t}$.
Hence $t^ngt^{-1}f = t^{n-1}f$.  It is now clear that (iii) holds.

(iv).  By Lemma~\ref{lem:really}\eqref{it:narb},  $(xyx)^{t^{n+2}} \in  \gen{x^4,\,y^4} = G_f = G_w$. 
\end{proof}

\begin{lemma}\label{lem:fixed} If Hypotheses~$\ref{hyp:example}$ hold, then  
$$V^{xyx} \quad =\quad \{t^nu \mid n \in \naturals -\{1\} \} 
\quad \cup \quad \{t^{n+2}w \mid n \in \naturals\}.$$
\end{lemma}

\begin{proof} Let $n \in \naturals$.

From~\cite[Definitions~I.3.4]{DicksDunwoody89}, we obtain the following.
\begin{align*}
&\iota^{-1}(t^nu) = t^nG_ue, 
&&\tau^{-1}(t^nu) = \emptyset,\\
&\iota^{-1}(t^{n+2}w) = t^{n+2}G_wf,
&& \tau^{-1}(t^{n+2}w) = t^{n+2}G_wt^{-2}e \,\,\cup \,\,t^{n+2}G_wt^{-1}f.
\end{align*}

By Lemma~\ref{lem:prefixed}(ii), (iii) and (iv),
the  edges of $T^{xyx}$ incident to $t^2w$    are $e$ and~$t^2f$,
the  edges of $T^{xyx}$ incident to $t^3w$  are $t^2f$ and~$t^3f$,
and, for $n \ge 2$, 
the  edges of $T^{xyx}$ incident to $t^{n+2}w$  are $t^ne$, $t^{n+1}f$ and $t^{n+2}f$.

Hence, in  $T^{xyx}$,
the  neighbours of $t^2w$  are $u$ and~$t^3w$,
the  neighbours of $t^3w$  are $t^2w$ and~$t^4w$,
and, for $n \ge 2$,
the  neighbours of $t^{n+2}w$  are $t^nu$, $t^{n+1}w$ and $t^{n+3}w$.

By Lemma~\ref{lem:prefixed}(i), if $n \ne 1$, then
the unique edge of $T^{xyx}$  incident to $t^{n}u$  is $t^{n}e$, 
and hence the unique neighbour of $t^{n}u$ in  $T^{xyx}$ is $t^{n+2}w$.

The result now follows.
\end{proof}

We now have the desired example.

\begin{theorem} There exists a group $G$ and 
a $G$-set $U$ such that $U$ is a $G$-retract of the vertex set of some $G$-tree 
but $U$ is not the vertex set of any $G$-tree.
\end{theorem}

\begin{proof} We assume that Hypotheses~$\ref{hyp:example}$ hold.

By Lemma~\ref{lem:bass-serre}, $U$ is a $G$-retract of the vertex set of some $G$-tree.

Suppose that there exists a $G$-tree  $T'$ with $VT' = U = Gu$.
We will derive a contradiction.  

Temporarily returning to the tree $T$, we let $L$ denote the subtree of $T$ 
with vertex set $\gen{t}w$ and edge set $\gen{t}f$.
Then $L$ is homeomorphic to $\reals$ and $t$ acts on $L$ by translation.
In particular, $\gen{t}$ acts freely on $VT$.  Hence,
$\gen{t}$ acts freely on $VT' \subseteq VT$.
As in~\cite[Proposition~I.4.11]{DicksDunwoody89},
there exists a subtree $L'$ of $T'$ homeomorphic to $\reals$
on which $t$ acts by translation.

 Let $v'$ denote the vertex of $L'$ closest to $u$ in~$T'$.
It is well known, and easy to prove,  that 
the $T'$-geodesic from $u$ to $t^2u$, denoted $T'[u,t^2u]$, 
 is the concatenation of the four $T'$-geodesics $T'[u,v']$, 
$T'[v', tv']$, $T'[tv',t^2v']$, and $T'[t^2v',t^2u]$.

By Lemma~\ref{lem:fixed}, and the fact that $\gen{t}$ acts freely  on $VT'$,
\begin{equation}\label{eq:fixed}
VT^{\prime xyx} = (Gu)^{xyx} = \{t^nu \mid n \in \naturals -\{1\}\}
= \{t^nu \mid n \in \naturals\} - \{tu\}.
\end{equation}

By~\eqref{eq:fixed}, or by direct calculation, $xyx$ fixes $u$, moves $tu$, and fixes
$t^2u$.  Thus,  $xyx$ fixes $T'[u,t^2u]$, 
and, hence, $xyx$ fixes $v'$, fixes $tv'$, and fixes~$t^2v'$. 

In particular, $tu \ne tv'$, hence $u \ne v'$, that is, $u \not\in L'$.  

Since $xyx$ fixes $v'$,  we see, by~\eqref{eq:fixed}, that $v' = t^nu$  
for some $n \in \naturals -\{1\}$.  Hence
$u = t^{-n}v' \in t^{-n}L' = L'.$ 
 This is a contradiction.  
\end{proof}

\medskip

\noindent{\textbf{\Large{Acknowledgments}}}

\medskip
\footnotesize

The research of the first-named author  was
funded by the DGI (Spain) through Project BFM2003-06613.

We are grateful to Gilbert Levitt for making us think about the 
sliding operation at a most opportune moment.  

We thank a referee
for several useful suggestions.

\vskip -.5cm\null

\bibliographystyle{amsplain}

\begin{thebibliography}{8}
\bibitem{delaHarpeValette89}
Pierre de la Harpe and Alain Valette, 
\newblock{\em La propri\'et\'e $(T)$ de Kazhdan pour les groupes localement compacts
 $($avec un appendice de Marc Burger$)$\/}, 
\newblock Ast\'erisque \textbf{175}, Soc.\ Math.\ de France, 
1989.
\vskip-0.44cm \null
\bibitem{Dicks94}
Warren Dicks,
\newblock {\em Equivalence of the strengthened Hanna Neumann conjecture and the 
amalgamated graph conjecture\/}, 
\newblock Invent.\ Math.\  \textbf{117}(1994), 373--389.
\newline Errata at http://mat.uab.cat/$\scriptstyle\sim$dicks/InvErr.html
\vskip-0.44cm \null
\bibitem{DicksDunwoody89}
Warren Dicks and M.~J. Dunwoody,
\newblock {\em Groups acting on graphs\/}, 
\newblock Cambridge Stud. Adv. Math.~\textbf{17}, CUP,  Cambridge, 1989.
\newline Errata at http://mat.uab.cat/$\scriptstyle\sim$dicks/DDerr.html
\vskip-0.44cm \null
\bibitem{DicksKropholler95}
Warren Dicks and Peter Kropholler,
\newblock {\em Free groups and almost equivariant maps\/}, 
\newblock {Bull.\ London Math.\ Soc.\ \textbf{27}(1995), 319--326.}
\newline  Addenda at http://mat.uab.cat/$\scriptstyle\sim$dicks/almost.html
\vskip-0.44cm \null
\bibitem{DicksLinnell06}
Warren Dicks and Peter A.\ Linnell,
\newblock {\em $L^2$-Betti numbers of one-relator groups\/}, 
\newblock{Math.\ Ann.\/} (to appear).
\newline  http://arxiv.org/abs/math.GR/0508370
\vskip-0.44cm \null
\bibitem{Dunwoody98}
M.\ J.\ Dunwoody,
\newblock {\em Folding sequences\/},  
\newblock{pp. 139--158 in: The Epstein birthday schrift
(eds. Igor Rivin, Colin Rourke and Caroline Series),
Geom.\ Topol.\ Monographs \textbf{1}\/},
\newblock{Geom.\ Topol.\ Publ.\/, Coventry, 1998.}
 \newline
http://www.maths.warwick.ac.uk/gt/GTMon1/paper7.abs.html
\vskip-0.44cm \null
\bibitem{Forester02}
Max Forester,
\newblock{\it Deformation and rigidity of simplicial
group actions on trees\/},
\newblock{Geom. Topol.\ \textbf{6}(2002), 219--267.} \newline
http://www.maths.warwick.ac.uk/gt/GTVol6/paper8.abs.html
\vskip-0.44cm \null
\bibitem{RipsSela97}
E.\ Rips and Z.\ Sela,
\newblock{\it Cyclic splittings of finitely presented groups and the canonical 
JSJ decomposition\/},
\newblock{Ann.\ Math.\ \textbf{146}(1997), 53--109.} 
\end{thebibliography}

\medskip

 \textsc{Warren Dicks,
Departament de  Matem\`atiques,
Universitat Aut\`onoma de Bar\-ce\-lo\-na,
E-08193 Bellaterra (Barcelona), Spain}

 \emph{E-mail address}{:\;\;}\texttt{dicks@mat.uab.cat}

\emph{URL}{:\;\;}\texttt{http://mat.uab.cat/$\scriptstyle\sim$dicks/}

\medskip

 \textsc{M.\ J.\ Dunwoody,
Department of Mathematics,
University of Southampton,
South\-amp\-ton, England  SO17 1BJ
}

 \emph{E-mail address}{:\;\;}\texttt{M.J.Dunwoody@maths.soton.ac.uk}

\emph{URL}{:\;\;}\texttt{http://www.maths.soton.ac.uk/staff/Dunwoody/}

\end{document}